\documentclass[11pt,reqno]{amsart}
\usepackage{indentfirst, amssymb,amsmath,amsthm}
\usepackage{enumitem}
\usepackage{xcolor}
\usepackage[colorlinks]{hyperref}
\usepackage{tikz-cd}
\newtheorem{theorem1}{Theorem}[section]
\newtheorem{example1}[theorem1]{Example}
\newtheorem{lemma1}[theorem1]{Lemma}

\newtheorem{remark1}[theorem1]{Remark}

\newtheorem{definition1}[theorem1]{Definition}

\newtheorem{note1}[theorem1]{Note}

\newcommand{\be}{\begin{equation}}
\newcommand{\ee}{\end{equation}}
\newcommand{\beas}{\begin{eqnarray*}}
\newcommand{\eeas}{\end{eqnarray*}}
\newcommand{\bea}{\begin{eqnarray}}
\newcommand{\eea}{\end{eqnarray}}

\numberwithin{equation}{section}

\begin{document}

\setlength{\unitlength}{1mm} \baselineskip .52cm
\setcounter{page}{1}
\pagenumbering{arabic}
\title[Lebesgue property and $G$-completeness in generalized quasi-uniform spaces]{Lebesgue property and $G$-completeness in generalized quasi-uniform spaces}

\author[Sugata Adhya and A. Deb Ray]{Sugata Adhya and A. Deb Ray}

\address{Department of Mathematics, The Bhawanipur Education Society College. 5, Lala Lajpat Rai Sarani, Kolkata 700020, West Bengal, India.}
\email {sugataadhya@yahoo.com}

\address{Department of Pure Mathematics, University of Calcutta. 35, Ballygunge Circular Road, Kolkata 700019, West Bengal, India.}
\email {debrayatasi@gmail.com}

\maketitle

\begin{abstract}
This paper extends the Lebesgue property and (weak) $G$-completeness to generalized quasi-uniform spaces. It investigates the connections between completeness, (weak) $G$-completeness, and the Lebesgue property of the product of generalized quasi-uniform spaces with those of the component spaces. It has been observed that the related behaviors differ for the Lebesgue property in contrast to completeness and (weak) $G$-completeness.
\end{abstract}


\section{\textbf{Introduction}}

The classes of metric spaces that lie between the classes of compact and complete metric spaces have fundamental importance in classical analysis. Such spaces help in well-understanding the ‘gap’ between compactness and completeness. During the last six decades those classes of metric spaces have been extensively studied and found applications in several branches of Mathematics (See e.g. \cite{l2} and references therein). Naturally, their study extended beyond the traditional metric framework. Lebesgue metric spaces and weak $G$-complete metric spaces are the examples of two such intermediate classes of metric spaces \cite{3greg2g-com, fm-uc, l2} that have been explored for fuzzy metric structure as well \cite{sal1,1sa-g-com,sal2,2greg-g-com}. This paper is devoted to extending such spaces in an uniform structure known as generalized quasi-uniform space \cite{rb}.

In point-set topology, uniform space found a natural place as a specialization of quasi-uniformity \cite{wilson}. Pervin established that not only a quasi-uniformity induces a topology but also each topological space finds a compatible quasi-uniform structure \cite{pervin1962quasi}. This result motivated Deb Ray and Bhowmik \cite{rb} to investigate a similar extended uniform structure for generalized topological spaces proposed by Császár \cite{cz}. In 2015, they extended the notion of quasi-uniformity to generalized quasi-uniformity so that each strong generalized topology can be made uniformizable in the new context \cite{rb}. 

Their work was extended in \cite{sa} with the study of product and completeness of such spaces. This paper serves as a subsequent continuation of that work. Here the study of Lebesgue property and (weak) $G$-completeness has been brought to the realm generalized quasi-uniform spaces.

After discussing necessary preliminaries in Section 2, here we extend Lebesgue property and (weak) \texorpdfstring{$G$}--completeness in generalized quasi-uniform spaces being motivated by their sequential descriptions through pseudo-Cauchy and $G$-Cauchy sequences respectively. We obtain varied forms of Lebesgue property when it is defined through nets and sequences. The natural question of finding the relationships between such variations of Lebesgue properties along with (weak) $G$-completeness and compactness has also been addressed. Subsequently we investigate the connections between completeness, (weak) $G$-completeness and Lebesgue properties of product of generalized quasi-uniform spaces with that of the component spaces. It has been observed that the related behaviour is different for Lebesgue property when compared to the cases of completeness and (weak) $G$-completeness.

\section{\textbf{Preliminaries}}

In this section, we recall some preliminaries on generalized topology and $g$-quasi uniformity which are to be required subsequently.

\begin{definition1}\normalfont\cite{c4e,cz,m3s}
Given a nonempty set $X,$ a subcollection $\mu$ of $\mathcal P(X)$ containning $\emptyset$ is called a generalized topology on $X$ if it is closed under arbitrary union. A generalized topology is said to be strong (or a supratopology) if it contains the entire set.
\end{definition1}

Elements of $\mu$ are called generalized open sets. A subset $A$ of $X$ is said to be a generalized neighborhood (or $\mu$-neighborhood) of $a$ if $a\in U\subset A$ for some $U\in\mu.$ For two generalized topological spaces $(X,\mu),~(Y,\mu'),$ a function $f: X\to Y$ is called generalized continuous if $f^{-1}(G)\in\mu,~\forall~G\in\mu'.$ 

\begin{definition1}
\normalfont\cite{wu-zhu} Let $\{(Z_i,\mu_i):i\in I\}$ be a family of generalized topological spaces. Then the generalized product topology\index{generalized!product topology defined by Wu and Zhu} on $Z=\prod\limits_{i\in I}Z_i$ is generated by $\{\pi_i:i\in I\},$ $\pi_i:Z\to Z_i,$ being the $i$-th projection mapping. In other words, it is the smallest generalized topology on $Z$ making each $\pi_i$ generalized continuous. Naturally, $\{\pi_i^{-1}(V_i):V_i\in\mu_i,i\in I\}$ forms a base for the generalized product topology on $Z.$
\end{definition1}

In 2015, Deb Ray and Bhowmik \cite{rb} extended the notion of quasi-uniformity so that each strong generalized topology can be made uniformizable in the new context.

\begin{definition1}
\normalfont\cite{rb} Given a set $X\neq\emptyset,~\mathcal U~(\neq\emptyset)\subset\mathcal P(X\times X)$ is called a generalized quasi-uniformity\index{generalized!quasi-uniformity} (or $g$-quasi uniformity\index{g@$g$-quasi!uniformity}) on $X$ and $(X,\mathcal U),$ a generalized quasi-uniform space\index{generalized!quasi-uniform space} (or $g$-quasi uniform space\index{g@$g$-quasi!uniform space}), if (i) $\Delta(X)\subset U,~\forall~U\in\mathcal U,$ (ii) $U\in\mathcal U$ and $V\supset U\implies V\in\mathcal U,$ (iii) $U\in\mathcal U\implies\exists~V\in\mathcal U$ such that $V\circ V\subset U.$ 

Given a $g$-quasi uniform space $(X,\mathcal U)$ and $\mathcal{B}~(\ne\emptyset)\subset\mathcal{U},$ $\mathcal{B}$ is called a base for $\mathcal U$ if given a member $V$ of $\mathcal U$ there is a member $B$ of $\mathcal B$ such that $B\subset V.$
\end{definition1}

\begin{theorem1}
\normalfont\cite{rb} Given a nonempty set $X$ and $\mathcal{B}~(\ne\emptyset)\subset\mathcal P(X\times X),$ $\mathcal{B}$ forms a base for some $g$-quasi uniformity on $X$ if and only if (i) $\Delta(X)\subset B,~\forall~B\in\mathcal B,$ and (ii) $B\in\mathcal B\implies\exists~V\in\mathcal B$ such that $V\circ V\subset B.$  

Moreover, such $\mathcal{B}$ is a base for the $g$-quasi uniformity $\{V\subset X:B\subset V\text{ for some }B\in\mathcal{B}\}$ on $X.$
\end{theorem1}

\begin{definition1}\normalfont\cite{rb}
For two $g$-quasi uniform spaces $(X,\mathcal U),~(Y,\mathcal U'),$ a mapping $f: X\to Y$ is called $g$-quasi uniform continuous\index{g@$g$-quasi!uniform continuous mapping} if $(f\times f)^{-1}(U)\in\mathcal U,~\forall~U\in\mathcal U'$.

It can be shown that $g$-quasi uniform continuity remains preserved under mapping composition.
\end{definition1}

\begin{theorem1}\label{th554}
\normalfont\cite{rb} For a $g$-quasi uniform space $(X,\mathcal U),~\mu(\mathcal U)=\{G\subset X:x\in G\implies U(x)\subset G\text{ for some }U\in\mathcal U\}$ forms a supratopology on $X.$ It is called the generalized\index{generalized!topology induced by a $g$-quasi uniformity} topology or supratopology generated (or induced) by $\mathcal U.$ Conversely, given a supratopology $\mu$ on a nonempty set $X,$ $\{(G\times G)\bigcup((X\backslash G)\times X):G\in\mu\}$ forms a base for some $g$-quasi uniformity\index{base for!a $g$-quasi uniformity} $\mathcal{U}_{\mu}$ on $X$ that generates $\mu$.
\end{theorem1}

\begin{theorem1}\label{thm554p}
\normalfont\cite{rb} Consider supratopological spaces $(X,\mu),$ $(X',\mu')$ and $g$-quasi uniform spaces $(Y,\mathcal U),$ $(Y',\mathcal U').$ Then 

(a) a mapping $f:(X,\mu)\to(X',\mu')$ is generalized continuous $\implies~f:(X,\mathcal U_\mu)\to(X',\mathcal U_{\mu'})$ is $g$-quasi uniform continuous,

(b) a mapping $g:(Y,\mathcal U)\to(Y',\mathcal U')$ is $g$-quasi uniform continuous $\implies~g:(Y,\mu(\mathcal U))\to(Y',\mu(\mathcal U'))$ is generalized continuous.
\end{theorem1}

\begin{theorem1}\label{th120223}\index{base for!product $g$-quasi uniformity}\index{base for!product $g$-quasi uniform space}
\normalfont\cite{sa} Let $\{(X_i,\mathcal U_i):i\in I\}$ be a family of $g$-quasi uniform spaces. Then $\mathcal B=\{(\pi_i\times\pi_i)^{-1}(U_i):U_i\in\mathcal U_i,i\in I\}$ forms a base for some $g$-quasi uniformity $\mathcal U$ on $X,$ where $X=\prod\limits_{i\in I}X_i.$
\end{theorem1}

\begin{definition1}
\normalfont\cite{sa} $\mathcal U$ (defined in Theorem \ref{th120223}) is called the product $g$-quasi uniformity\index{product!g@$g$-quasi uniformity} (to be denoted by $\prod\limits_{i\in I}\mathcal U_i)$ on $X$ and the pair $(X,\mathcal U),$ the product $g$-quasi uniform space\index{product!g@$g$-quasi uniform space}, for the family $\{(X_i,\mathcal U_i):i\in I\}.$
\end{definition1}

We finish this section by recalling certain preliminaries on Lebesgue property and weak $G$-completeness for metric spaces.

\begin{definition1}
\normalfont\cite{l2} A metric space on which every real-valued continuous function is uniformly continuous is said to be Lebesgue (or Atsuji space). 
\end{definition1}

\begin{definition1}
\normalfont\cite{l3} A sequence $(x_n)$ in a metric space $(X,d)$ is said to be pseudo-Cauchy if given $\epsilon>0,k\in\mathbb N$ there exists distinct $m,n~(>k)\in\mathbb N$ such that $d(x_m,x_n)<\epsilon.$
\end{definition1}

\begin{theorem1}
\normalfont\cite{l2,l3} A metric space is Lebesgue if and only if every pseudo-Cauchy sequence having distinct terms clusters in it. 
\end{theorem1}

\begin{definition1}
\normalfont\cite{1sa-g-com,2greg-g-com,3greg2g-com} A sequence $(x_n)$ in a metric space $(X,d)$ is called $G$-Cauchy\index{g@$G$-Cauchy sequence!in metric spaces} if $\lim\limits_{n\to\infty}d(x_{n+p},x_n)=0,~\forall~ p\in\mathbb N$ (or equivalently, $\lim\limits_{n\to\infty}d(x_{n+1},x_n)=0$). A metric space in which every $G$-Cauchy sequence converges is said to be weak $G$-complete.
\end{definition1}

Both Lebesgue property and weak $G$-completeness are strictly intermediate between compactness and completeness of metric spaces \cite{2greg-g-com,l2}.

\section{\textbf{Lebesgue Property and (Weak) \texorpdfstring{$G$}--Completeness in \texorpdfstring{$g$}--Quasi Uniform Spaces}\label{gquasi-defi-Lebesgue}}

In this section, we explore the interrelation among varied forms of Lebesgue property and (weak) \texorpdfstring{$G$}--completeness in \texorpdfstring{$g$}--quasi uniform spaces. We begin this section by introducing ideas related to net convergence in generalized topology and generalized quasi-uniformity.

\begin{definition1}
\normalfont\cite{kelly} Let $Z$ be a nonempty set. A mapping $S$ from an ordered set $(D,\ge)$ to $Z$ is called a net\index{net} in $Z$. For convenience we denote it by $(S_n)_{n\in D},$ where $S_n$ signifies the image of $n$ under $S.$

Unless stated otherwise, $D$ denotes an ordered set $(D,\ge)$ throughout the chapter.
\end{definition1}

\begin{definition1}\label{def558}
\normalfont Let $(X,\mu)$ be a generalized topological space and $(S_n)_{n\in D},$ a net in it. Then

(a) $(S_n)_{n\in D}$ is said to be convergent\index{convergent net} to $c\in X$ if given a $\mu$-neighborhood $U$ of $c,$ there exists $m\in D$ such that $S_n\in U,~\forall~n\ge m;$ 

(b) a point $c\in X$ is called a cluster point\index{cluster point} of $(S_n)_{n\in D}$ if given a $\mu$-neighborhood $U$ of $c$ and $n\in D,$ there exists $m\in D$ with $m\ge n$ such that $S_m\in U.$
\end{definition1}

\begin{definition1}
\normalfont Let $(S_n)_{n\in D}$ be a net in a $g$-quasi uniform space $(Z,\mathcal U)$ and $c\in Z.$ Then

(i) $c$ is called a cluster point\index{cluster point} of $(S_n)_{n\in D}$ in $(Z,\mathcal U)$ if it is so in $(Z,\mu(\mathcal U));$

(ii) $(S_n)_{n\in D}$ is said to be convergent\index{convergent net} to $c$ in $(Z,\mathcal U)$ if it is so in $(Z,\mu(\mathcal U)).$

Clearly if $(S_n)_{n\in D}$ is convergent to $c,$ then $c$ is a cluster point of $(S_n)_{n\in D}.$
\end{definition1}

\begin{definition1}
\normalfont Let $(Z,\mathcal U)$ be a $g$-quasi uniform space.

(i) A net $(S_n)_{n\in D}$ in $Z$ is called Cauchy \cite{kelly}\index{Cauchy net} if for $U\in\mathcal U$ there exists $m\in D$ such that $(S_p,S_q)\in U,~\forall~p,q\ge m;$  

(ii) A net $(S_n)_{n\in D}$ in $Z$ is called pseudo-Cauchy\index{pseudo-Cauchy!net in $g$-quasi uniform spaces} if for $U\in\mathcal U,$ $p\in D$ there exist $m,n~(m\ne n)\in D$ with $m,n\ge p$ such that $(S_m,S_n)\in U;$

(iii) A sequence $(x_n)$ in $Z$ is called $G$-Cauchy\index{g@$G$-Cauchy sequence! in $g$-quasi uniform spaces} if for $U\in\mathcal U$ there exists $k\in\mathbb N$ such that $(x_n,x_{n+1})\in U,~\forall~n\ge k.$
\end{definition1}

\begin{definition1}
\normalfont A $g$-quasi uniform space $(Z,\mathcal U)$ is called

(i) complete \cite{kelly}\index{complete!g@$g$-quasi uniform space} if every Cauchy net converges to some point in it;

(ii) Lebesgue\index{Lebesgue!g@$g$-quasi uniform space} if every pseudo-Cauchy net $(S_n)_{n\in D}$ having distinct terms (i.e., $S_m\ne S_n,~\forall~m\ne n)$ has a cluster point in it;

(iii) strongly Lebesgue\index{strongly Lebesgue!g@$g$-quasi uniform space} if every pseudo-Cauchy net has a cluster point in it;

(iv) sequentially Lebesgue\index{sequentially!Lebesgue $g$-quasi uniform space} if every pseudo-Cauchy sequence having distinct terms has a cluster point in it;

(v) $G$-complete\index{g@$G$-complete!g@$g$-quasi uniform space} if every $G$-Cauchy sequence is convergent to some point in it;

(vi) weak $G$-complete\index{weak $G$-complete!g@$g$-quasi uniform space} if every $G$-Cauchy sequence has a cluster point in it.
\end{definition1}

\begin{definition1}
\normalfont A $g$-quasi uniform space $(Z,\mathcal U)$ is said to be compact\index{compact!g@$g$-quasi uniform space} if $(Z,\mu(\mathcal U))$ is compact.
\end{definition1}

Following \cite{folland}, we can show that every net in a compact $g$-quasi uniform space has a cluster point. Consequently, we have the following:

\begin{theorem1}
\normalfont A compact $g$-quasi uniform space is strongly Lebesgue.
\end{theorem1}

In view of the foregoing discussion, we have the following chain of implications for $g$-quasi uniform spaces:\\

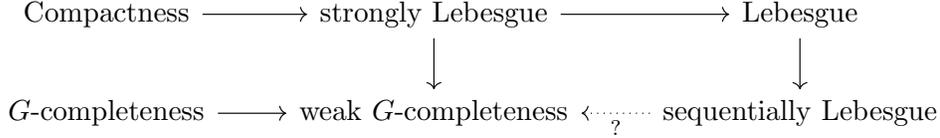
\noindent\begin{tikzcd}
\text{Compactness} \arrow{r} & \text{strongly Lebesgue} \arrow{r} \arrow{d} & \text{Lebesgue} \arrow{d} \\
    G\text{-completeness} \arrow{r} & \text{weak }G\text{-completeness} & \arrow[l,dotted,"?"] \text{sequentially Lebesgue}
\end{tikzcd}\\

\begin{note1}
\normalfont A Lebesgue (resp. sequentially Lebesgue) $g$-quasi uniform space may not be strongly Lebesgue. For example, consider the $g$-quasi uniform space $(\mathbb Z,\mathcal{U})$ where $\mathcal{U}$ is the uniformity (and hence a $g$-quasi uniformity) on $\mathbb Z$ generated by the discrete metric on it. Then $(\mathbb Z,\mathcal{U})$ is Lebesgue (resp. sequentially Lebesgue). However it is not strongly Lebesgue since $\{1,1,2,2,3,3,\cdots\}$ is a pseudo-Cauchy sequence in it without any cluster point.
\end{note1}

It is clear from the last diagram that a condition that makes sequentially Lebesgue $g$-quasi uniform spaces weak $G$-complete is desirable in the present context. The following theorem discusses such a condition.

\begin{theorem1}\label{th514countable}
\normalfont Let $(X,\mathcal{U})$ be a sequentially Lebesgue $g$-quasi uniform space. If $\mathcal{U}$ has a countable base, then $(X,\mathcal{U})$ is weak $G$-complete.
\end{theorem1}

\begin{proof}
\normalfont Let $\mathcal{B}=\{B_n:n\in\mathbb N\}$ be a countable base for $\mathcal{U}.$ Choose a $G$-Cauchy sequence $(x_n)$ in $(X,\mathcal{U}).$ 

If $(x_n)$ has a constant subsequence, then it clearly has a cluster point in $(X,\mathcal{U}).$ So we assume that $(x_n)$ has no constant subsequence. 

Consider the infinite matrix $(X_{ij})$ where $X_{ij} =
\begin{cases}
B_j,  & \text{if $i\ge j$} \\
\emptyset, & \text{otherwise}
\end{cases}.$

That is,

$$(X_{ij})=\begin{pmatrix}
B_1 & \emptyset & \emptyset & \emptyset & \cdots \\
B_1 & B_2 & \emptyset & \emptyset & \cdots\\
B_1 & B_2 & B_3 & \emptyset & \cdots\\
\vdots & \vdots & \vdots & \vdots & \ddots & 
\end{pmatrix}.$$

We construct a pseudo-Cauchy subsequence of $(x_n)$ by induction.

Choose $r_{111}\in\mathbb N$ such that $x_{r_{111}},x_{r_{111}+1}\in X_{11}$ and $x_{r_{111}}\ne x_{r_{111}+1}.$ Set $x_{r_{112}}=x_{r_{111}+1}$

In view of the fact that $(x_n)$ has no constant subsequence, for the choices of an $m\in\mathbb N$ and $x_{r_{ij1}}, x_{r_{ij2}}\in X_{ij},$ $\forall~j\le i\le m,$ we choose $x_{r_{(m+1)j1}}, x_{r_{(m+1)j2}}\in X_{(m+1)j},$ $\forall~j\le m+1$ such that

(i) $ijk<qst$ (in dictionary order) $\implies$ 
$r_{ijk}<r_{qst},~\forall~i,q=1,2,\cdots,m+1;$ $j\le i,s\le q$ and $k,t=1,2;$

(ii) $x_{r_{ijk}}$’s are distinct, $\forall~i=1,2,\cdots,m+1;$ $j\le i $ and $k =1,2.$

It is clear that the dictionary order of $ijk$ makes $(r_{ijk})$ a strictly increasing sequence of natural numbers. By renaming the sequence as $(r_n),$ we obtain a subsequence $(x_{r_n})$ of $(x_n)$ having distinct terms.

Choose $V\in\mathcal{U}$ and $p\in\mathbb N.$ Since $\mathcal{B}$ forms a base for $\mathcal{U},$ there exists $q\in\mathbb N$ such that $B_q\subset V.$ It is clear from the construction of $(x_{r_n})$ that there exists $m\in\mathbb N$ with $m>p$ such that $(x_{r_m},x_{r_{m+1}})\in B_q\subset V.$ Since $r_m,r_{m+1}>p,$ $(x_{r_n})$ is pseudo-Cauchy and hence, it has a cluster point $c$ in $(X,\mathcal{U}).$ Consequently $c$ is a cluster point of $(x_n)$ in $(X,\mathcal{U}).$ 

Thus $(X,\mathcal{U})$ is weak $G$-complete.\end{proof}

Thus for $g$-quasi uniform spaces having countable bases, we can improve the previous chains of implications as follows:\\

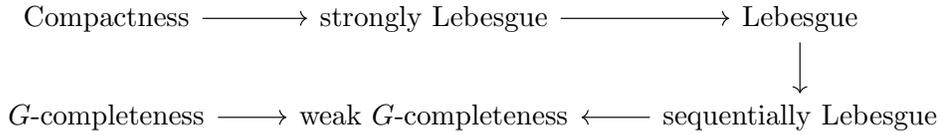
\noindent\begin{tikzcd}
\text{Compactness} \arrow{r} & \text{strongly Lebesgue} \arrow{r} & \text{Lebesgue} \arrow{d} \\
    G\text{-completeness} \arrow{r} & \text{weak }G\text{-completeness} & \arrow{l} \text{sequentially Lebesgue}
\end{tikzcd}\\

\section{\textbf{Lebesgue Property and (Weak) $G$-Completeness in Product of $g$-Quasi Uniform Spaces}}

We begin by recalling the following results from \cite{sa}:

\begin{lemma1}\label{130223}
\normalfont Let $(X,\mathcal U)$ be the product for a family $\{(X_i,\mathcal U_i):i\in I\}$ of $g$-quasi uniform spaces. Then a net $(R_n)_{n\in D}$ in $X$ is Cauchy if and only if for each $i\in I,$ the net $(\pi_i(R_n))_{n\in D}$ in $X_i$ is Cauchy.
\end{lemma1}

\begin{theorem1}\label{chap5_complete_thm}
\normalfont Let $(X,\mathcal U)$ be the product for a family $\{(X_i,\mathcal U_i):i\in I\}$ of $g$-quasi uniform spaces. Then $(X,\mathcal U)$ is complete if and only if $(X_i,\mathcal U_i)$ is complete, $\forall~i\in I.$
\end{theorem1}

In what follows, we examine similar cases for varied forms of Lebesgue and $G$-complete $g$-quasi uniform spaces.

\begin{lemma1}\label{chap5_psedo-Cauchy_lemma}
\normalfont Let $(X,\mathcal U)$ be the product for a family $\{(X_i,\mathcal U_i):i\in I\}$ of $g$-quasi uniform spaces. Then a net $(R_n)_{n\in D}$ in $X$ is pseudo-Cauchy if and only if for each $i\in I,$ the net $(\pi_i(R_n))_{n\in D}$ in $X_i$ is pseudo-Cauchy.
\end{lemma1}

\begin{proof}
\normalfont Let $(R_n)_{n\in D}$ be a pseudo-Cauchy net in $X$. Then for chosen $i\in I, U_i\in\mathcal U_i$ and $q\in D$ there is $m,n~(m\ne n)\in D$ with $m,n\ge q$ such that $(R_m,R_n)\in(\pi_i\times\pi_i)^{-1}(U_i).$ i.e., $(\pi_i(R_m),\pi_i(R_n))\in U_i\implies(\pi_i(R_n))_{n\in D}$ is pseudo-Cauchy in $X_i.$

\textit{Conversely}, suppose $(R_n)_{n\in D}$ be a net in $X$ such that $(\pi_i(R_n))_{n\in D}$ pseudo-Cauchy in $X_i,$ for each $i\in I$. 

Choose $U\in\mathcal U,p\in D.$ Then $(\pi_i\times\pi_i)^{-1}(U_i)\subset U$ for some $i\in I, U_i\in\mathcal U_i.$ 

Due to the hypothesis, there exist $m,n~(m\ne n)\in D$ with $m,n\ge p$ such that $(\pi_i(R_m),\pi_i(R_n))\in U_i.$ That is $(R_m,R_n)\in(\pi_i\times\pi_i)^{-1}(U_i)\subset U.$ 

Thus $(R_n)_{n\in D}$ is pseudo-Cauchy in $X.$
\end{proof}

Similar chain of arguments yield identical results for $G$-Cauchy sequences. We state the result without proof.

\begin{lemma1}
\normalfont Let $(X,\mathcal U)$ be the product for a family $\{(X_i,\mathcal U_i):i\in I\}$ of $g$-quasi uniform spaces. Then a sequence $(x_n)$ in $X$ is $G$-Cauchy if and only if for each $i\in I,$ the sequence $(\pi_i(x_n))$ in $X_i$ is $G$-Cauchy.
\end{lemma1}

\begin{theorem1}
\normalfont Let $(X,\mathcal U)$ be the product for a family $\{(X_i,\mathcal U_i):i\in I\}$ of $g$-quasi uniform spaces. Then $(X,\mathcal U)$ is strongly Lebesgue if and only if $(X_i,\mathcal U_i)$ is strongly Lebesgue, $\forall~i\in I.$
\end{theorem1}

\begin{proof}
\normalfont Let $(X,\mathcal U)$ be strongly Lebesgue. Choose a pseudo-Cauchy net $(S_n)_{n\in D}$ in $X_i$ where $i\in I.$ 
Fix $c_k\in X_k,~\forall~k\in I\backslash\{i\}.$ Then the net $(R_n)_{n\in D}$ in $X,$ defined as in Theorem \ref{chap5_complete_thm}, is pseudo-Cauchy by Lemma \ref{chap5_psedo-Cauchy_lemma}. Consequently $(R_n)_{n\in D}$ has a cluster point $x$ in $X.$

Choose a generalized neighborhood $G_i$ of $\pi_i(x)$ in $X_i.$ Since $\pi_i$ is generalized continuous, given $q\in D$ there exists $k\in D$ with $k\ge q$ such that $R_k\in\pi_i^{-1}(G_i)$ and hence, $\pi_i(R_k)\in G_i.$ Thus $\pi_i(x)$ is a cluster point of $(S_n)_{n\in D}=(\pi_i(R_n))_{n\in D}$ in $X_i.$ Hence $(X_i,\mathcal U_i)$ is strongly Lebesgue.

\textit{Conversely,} let $(X_i,\mathcal U_i)$ be strongly Lebesgue, $\forall~i\in I.$ Choose a pseudo-Cauchy net $(R_n)_{n\in D}$ in $X.$ Then for each $i\in I,(\pi_i(R_n))_{n\in D}$ is pseudo-Cauchy in $X_i$ by Lemma \ref{chap5_psedo-Cauchy_lemma}.

Due to the hypothesis, for chosen $i\in I$ there exists $x_i\in X_i$ such that $x_i$ is a cluster point of $(\pi_i(R_n))_{n\in D}$ in $X_i.$ Choose a generalized neighborhood $G$ of $x=(x_i)_{i\in I}$ in $X$ and $q\in D.$ Then $V(x)\subset G$ for some $V\in\mathcal U.$ So there exist $j\in I,V_j\in\mathcal U_j$ such that $((\pi_j\times\pi_j)^{-1}(V_j))(x)\subset G.$ That is, $\pi_j^{-1}(V_j(x_j))\subset G.$ Since $V_j(x_j)$ is a generalized neighborhood of $x_j,$ there exists $k\in D$ with $k\ge q$ such that $\pi_j(R_k)\in V_j(x_j)$ and hence, $R_k\in G.$ 

So $x=(x_i)_{i\in I}$ is a cluster point of $(R_n)_{n\in D}$ in $X.$ 

Thus $(X,\mathcal U)$ is strongly Lebesgue.
\end{proof}

\begin{theorem1}\label{leb_prod}
\normalfont Let $(X,\mathcal U)$ be the product for a family $\{(X_i,\mathcal U_i):i\in I\}$ of $g$-quasi uniform spaces. If $(X,\mathcal U)$ is Lebesgue then so is $(X_i,\mathcal U_i),$ for each $i\in I.$
\end{theorem1}

\begin{proof}
\normalfont Let $(X,\mathcal U)$ be Lebesgue. Choose $i\in I$ and a pseudo-Cauchy net $(S_n)_{n\in D}$ having distinct terms in $X_i.$

Fix $c_k\in X_k,$ $\forall~k\in I\backslash\{i\}.$ Then the net $(R_n)_{n\in D}$ in $X$ (defined as before) is pseudo-Cauchy, by Lemma \ref{chap5_psedo-Cauchy_lemma}, having distinct terms. Consequently, $(R_n)_{n\in D}$ has a cluster point $x$ in $X.$ Hence, proceeding as before, we see that $\pi_i(x)$ is a cluster point of $(S_n)_{n\in D}$ in $X_i.$ Thus $(X_i,\mathcal U_i)$ is Lebesgue.
\end{proof}

Proceeding similarly as before we can prove the following theorem:

\begin{theorem1}\label{seq_leb_prod}
\normalfont Let $(X,\mathcal U)$ be the product for a family $\{(X_i,\mathcal U_i):i\in I\}$ of $g$-quasi uniform spaces. If $(X,\mathcal U)$ is sequentially Lebesgue then so is $(X_i,\mathcal U_i),$ for each $i\in I.$
\end{theorem1}

The converses of Theorem \ref{leb_prod} and Theorem \ref{seq_leb_prod} are not true, in general. In support, we provide the following example.

\begin{example1}\label{leb_prod_exm}
\normalfont Consider the collection $\{(X_i,\mathcal{U}_i):i\in\mathbb N\}$ of $g$-quasi uniform spaces where $\forall~i\in\mathbb N,$ $X_i=\mathbb N\bigcup\left\{\frac{1}{i}\right\}$ and $\mathcal{U}_i$ is the uniformity (and hence a $g$-quasi uniformity) on $X_i$ induced by the usual metric on it. Clearly $\triangle(X_i)\in\mathcal{U}_i,~\forall~i\in\mathbb N.$

We first show that $(X_i,\mathcal{U}_i)$ is Lebesgue (and hence sequentially Lebesgue), $\forall~i\in\mathbb N.$ Choose $i\in\mathbb N.$ If possible, let $(R_n)_{n\in D}$ be a pseudo-Cauchy net in $X_i$ having distinct terms. Since $\triangle(X_i)\in\mathcal{U}_i,$ for given $p\in D$ there exist $m,n~(m\ne n)\in D$ with $m,n\ge p$ such that $(R_m,R_n)\in\triangle(X_i)\implies R_m=R_n,$ a contradiction. Thus there is no pseudo-Cauchy net in $X_i$ having distinct terms. Hence $(X_i,\mathcal{U}_i)$ is Lebesgue (and hence sequentially Lebesgue), $\forall~i\in\mathbb N.$

Let $(X,\mathcal{U})$ be the product $g$-quasi uniform space for the family $\{(X_i,\mathcal{U}_i):i\in\mathbb N\}.$ We show that $(X,\mathcal{U})$ is not sequentially Lebesgue (and hence not Lebesgue).

Consider the sequence $(x_n)$ in $X$ where\\ 
    $x_1=(0,0,0,0,0,\cdots),\\
    x_2=(1,0,0,0,0,\cdots),\\
    x_3=(2,0,0,0,0,\cdots),\\
    x_4=(2,\frac{1}{2},0,0,0,\cdots),\\
    x_5=(3,0,0,0,0,\cdots),\\
    x_6=(3,0,\frac{1}{3},0,0,\cdots),\\
    x_7=(4,0,0,0,0,\cdots),\\
    x_8=(4,0,0,\frac{1}{4},0,\cdots),\\$
and so on. Clearly $(x_n)$ is a pseudo-Cauchy sequence in $X$ having distinct terms. 

However $(x_n)$ has no cluster point in $X.$ If not, let $y=(y_1,y_2,\cdots)$ be  a cluster point of $(x_n)$ in $X.$ Then there exists $k\in\mathbb N$ such that $\pi_1(x_k)>y_1.$ Consequently $\pi_1(x_n)>y_1,~\forall~n\ge k.$ Since $((\pi_1\times\pi_1)^{-1}(\triangle(X_1))(y)$ is a neighborhood of $y$ in $X$ there exists $q\in\mathbb N$ with $q\ge k$ such that $x_q\in((\pi_1\times\pi_1)^{-1}(\triangle(X_1))(y).$ Thus $\pi_1(x_q)=y_1,$ a contradiction. 

Hence $(X,\mathcal{U})$ is not sequentially Lebesgue (and hence not Lebesgue).
\end{example1}

Following the chain of arguments used to prove earlier results, we find that the $G$-completeness and weak $G$-completeness of the product $g$-quasi uniform space is achieved via those of the component spaces and conversely. We state the results here without proof:

\begin{theorem1}
\normalfont Let $(X,\mathcal U)$ be the product for a family $\{(X_i,\mathcal U_i):i\in I\}$ of $g$-quasi uniform spaces. Then $(X,\mathcal U)$ is $G$-complete if and only if $(X_i,\mathcal U_i)$ is $G$-complete, for each $i\in I.$
\end{theorem1}

\begin{theorem1}\label{weak_g}
\normalfont Let $(X,\mathcal U)$ be the product for a family $\{(X_i,\mathcal U_i):i\in I\}$ of $g$-quasi uniform spaces. Then $(X,\mathcal U)$ is weak $G$-complete if and only if $(X_i,\mathcal U_i)$ is weak $G$-complete, for each $i\in I.$
\end{theorem1}

\begin{remark1}
\normalfont We see that each $(X_i,\mathcal{U}_i)$ in Example \ref{leb_prod_exm} is Lebesgue and hence, weak $G$-complete by Theorem \ref{th514countable}. Thus, by Theorem \ref{weak_g}, their product $(X,\mathcal{U})$ is weak $G$-complete. However $(X,\mathcal{U})$ is not sequentially Lebesgue. Thus a weak $G$-complete $g$-quasi uniform space may not be (sequentially) Lebesgue.
\end{remark1}

\end{document}